\documentclass[a4paper, 12pt]{article}
\usepackage{amsmath}
\usepackage{amssymb}
\usepackage{psfrag}
\usepackage{latexsym}
\usepackage{enumerate}
\usepackage{graphics}
\usepackage{subfigure}
\usepackage{caption}
\input epsf
\usepackage[dvips]{graphicx}
\usepackage{enumerate}
\usepackage[english]{babel}
\usepackage{amsmath,amssymb,amsthm}

\newcommand{\be}{\begin{equation}}
\newcommand{\ee}{\end{equation}}
\newcommand{\bea}{\begin{eqnarray}}
\newcommand{\eea}{\end{eqnarray}}

\newcommand{\nod}{\noindent}

\newcommand{\ba}{\begin{array}}
\newcommand{\ea}{\end{array}}
\newcommand{\bc}{\begin{center}}
\newcommand{\ec}{\end{center}}

\allowdisplaybreaks

\begin{document}
\title{\bf Exact and approximate epidemic models on networks: a new, improved closure relation\\}

\author{Istvan Z. Kiss $^{1,\ast}$ \& P\'eter L. Simon $^{2}$}

\maketitle

\begin{center}

$^1$ School of Mathematical and Physical Sciences, Department of
Mathematics, University of Sussex, Falmer,
Brighton BN1 9QH, UK\\
$^2$ Institute of Mathematics, E\"otv\"os Lor\'and University
Budapest, Budapest, Hungary\\

\end{center}

\vspace{8cm}
\begin{flushleft}
$\ast$ corresponding author\\
email: i.z.kiss@sussex.ac.uk\\
\end{flushleft}

\newpage
\begin{abstract}
Recently, research that focuses on the rigorous understanding of the relation between simulation
and/or exact models on graphs and approximate counterparts has gained lots of momentum. This
includes revisiting the performance of classic pairwise models with closures at the level of pairs
and/or triples as well as effective-degree-type models and those based on the probability
generating function formalism. In this paper, for a fully connected graph and the simple $SIS$
(susceptible-infected-susceptible) epidemic model, a novel closure is introduced. This is done via
using the equations for the moments of the distribution describing the number of infecteds at all
times combined with the empirical observations that this is well described/approximated by a
binomial distribution with time dependent parameters. This assumption allows us to express higher
order moments in terms of lower order ones and this leads to a new closure. The significant
feature of the new closure is that the difference of the exact system, given by the
Kolmogorov equations, from the solution of the newly defined approximate system is of order
$1/N^2$. This is in contrast with the $\mathcal{O}(1/N)$  difference corresponding to the
approximate system obtained via the classic triple closure.\\
\end{abstract}

\nod {\bf Keywords:} Markov Chain; epidemic; pairwise model;
pairwise closure

\newpage

\section{Introduction: exact stochastic models on and\\ of graphs}

In this Section we present two important examples that motivate our investigations and a possible
extension where progress can be made following some of the methods introduced in this paper. Here,
we use a dynamical system type approach, where the Kolmogorov equations are simply considered as a
system of linear ODEs with a transition rate matrix with specific properties such as special
tri-diagonal structure and/or well defined functional form for the transmission rates. For
example, consider a Markov chain with finite state space $\{ 0, 1, \ldots ,N\}$ and denote by
$p_k(t)$ the probability that the system is in state $k$ at time $t$ (with a given initial state
that is not specified at the moment). Assuming that starting from state $k$ the system can move to
either state $k-1$ or to state $k+1$, the Kolmogorov equations of the Markov chain take the form
\begin{equation}
    \dot{p}_k=a_{k-1}p_{k-1}-(a_{k}+c_{k})p_k+c_{k+1}p_{k+1},\quad k=0,\ldots,
    N.
    \tag{KE} \label{eq:Kolmogorov}
\end{equation}

The first motivation for our study comes from epidemiology where a paradigm disease transmission
model is the simple susceptible-infected -susceptible ($SIS$) model on a completely
connected graph with $N$ nodes, i.e. all individuals are connected to each other. From the disease
dynamic viewpoint, each individual is either susceptible ($S$) or infected ($I$) -- the
susceptible ones can be infected at rate $\tau$ if connected to an infected node, and the infected
ones can recover at rate $\gamma$ and become susceptible again. It is known that in this case the
$2^N$-dimensional system of Kolmogorov equations can be lumped to a $N+1$-dimensional system, see
\cite{ExactEpiMod}.

The lumped Kolmogorov equations take again the form
(\ref{eq:Kolmogorov}) with
\begin{equation}
    a_{k}=\tau k(N-k), \quad c_{k}=\gamma k \quad \text{for} \quad k=0,\ldots, N \quad \text{with} \quad
    a_{-1}=c_{N+1}=0.
   \label{eq:SIS_coeff_fully_con_graph}
\end{equation}
For homogeneous random graphs, where every node has $n(<N)$ links to
other nodes in the network, a similar system can be written down.
In this case the transition rates are as follows,
\begin{equation}
    a_{k}=\tau nk\frac{N-k}{N-1}, \quad c_{k}=\gamma k \quad \text{for} \quad k=0,\ldots, N \quad \text{with} \quad
    a_{-1}=c_{N+1}=0.
    \label{eq:SIS_coeff_hom_ran_graph}
\end{equation}
This equation is not exact in the sense that it is only an
approximation to the exact process unfolding on a homogeneous random
graph. This is not of immediate relevance as in this paper we
concentrate on approximations of the Kolmogorov equations via
low-dimensional ODEs, in particular we aim to derive approximate
equations for the moments of the stochastic process described by
the Kolmogorov equations.

A final example stems from a recent model of a dynamic network
where links can be activated and deleted at random while subjected
to a global constraint on the total number of links in the
network (see Kiss et al. \cite{KissDynNetw}). This example
is a prerequisite to models where dynamics
on the network is coupled with the dynamic of the network Gross and Blasius \cite{ThiloAdCoevNetw}.
In the case of epidemic propagation on networks it is straightforward to assume that the
propagation of the epidemic has an effect on the structure of the
network. For example, susceptible individuals try to cut their
links in order to minimise their exposure to infection. This leads
to a change in network structure which in turn impacts on how the
epidemic spreads. The first step in modelling this phenomenon is
an appropriate dynamic network model such as the recently proposed
globally-constrained Random Link Activation-Deletion (RLAD) model.
This can be described in terms of Kolmogorov equations as follows,
\begin{eqnarray*}
\dot p_k(t)&=& \alpha (N-(k-1))\left(1-\frac{k-1}{K_1^{max}}\right) p_{k-1}(t)\\
&-& \left[\alpha (N-k)\left(1-\frac{k}{K_1^{max}}\right)+\omega k \right] p_k(t)+\omega (k+1) p_{k+1}(t),\\
k&=&0,\ldots, N,
\label{eq:RLAD}
\end{eqnarray*}
where $p_k(t)$ denotes the probability that at time $t$ there are
$k$ active links in the network with $N$ denoting the total number of
potential edges. It is assumed that non-active links are activated
independently at random at rate $\alpha$ and that existing links
are broken independently at random at rate $\omega$. Furthermore,
the link creation is globally constrained by introducing a
carrying capacity $K_{1}^{max}$, that is the network can only
support a certain number of edges as given by $K_{1}^{max}$.

Using the above notation, here
\begin{equation}
    a_{k}=\alpha(N-k)\left(1-\frac{k}{K_1^{\max}}\right), \quad c_{k}=\omega k, \quad \text{for} \quad k=0,\ldots, N \quad \text{with} \quad
    a_{-1}=c_{N+1}=0.
    \label{eq:rewRLAD_coeff}
\end{equation}
The common ingredient of all the models above is the set of Kolmogorov equations.
Solving these even numerically, more often than not, is challenging
and impossible simply due to the large number of equations. This has led to various approaches that
set out to derive approximate models that are low-dimensional and can capture the
exact dynamics in terms of the expected values of some well-defined quantities.
These approaches could be broadly classified into heuristic and rigorous or semi-rigorous.
For example, the original derivation of pairwise models was based on heuristic arguments \cite{Keeling1999}.
More rigorous approaches include the lumping technique that exploits system symmetries and allows a significant reduction
of the  state space \cite{ExactEpiMod}. Additional models such as the effective degree models \cite{BallNeal, EffDegMod} and the probability
generating function approach by Volz \cite{Miller, Volz} led to excellent agreement with simulations
without the explicit proof of convergence but with strong
probabilistic arguments based on results by Kurtz \cite{Kurtz70}.

In this paper, starting from the Kolmogorov equations, given by the simple $SIS$ epidemic model on a fully connected graph, the evolution equations for the moments are derived and are interpreted in terms of and compared to the classic pairwise equations.
The equations for the moments are not self-contained, and a new/novel closure is proposed. This is based on the fact that the number of infecteds at all times is well approximated by a binomial distribution with time dependent parameters. The performance of the new closure is investigated numerically in the closing Section and we show that this is superior to the classic closure at the level of triples used for pairwise models.

\section{Pairwise models: closure relations and their performance}

It is of great interest to understand how and when the full set of Kolmogorov equations can be
approximated via low-dimensional ODEs and also to assess the performance of the approximate models
by, for example, working out the rate of convergence of the exact/full model towards the solution
of the approximate or mean-field model or simply estimating the absolute difference between the two.
Pairwise models have been widely used as an alternative or
as a companion to simulations of epidemic models on mainly homogeneous random graphs. While
originally, the pairwise equations have been heuristically defined, more recently, Simon et al.
\cite{ExactEpiMod} and Taylor et al. \cite{Taylor} have shown that these can be derived directly
from the exact Kolmogorov equations and that these pairwise equations are exact before closure on
arbitrary graphs. Focusing on the simple $SIS$ type model, the first moment of the distribution is
given by
\begin{equation}
[\dot I]=\tau[SI]-\gamma[I],
\end{equation}
where $[I](t)=\sum_{k=0}^{N}kp_{k}(t)$. This is not a closed
equation since $[SI]$ itself is a variable and an equation for this
is needed. However, we can look for an approximation whereby the
expected number of edges $[SI] \simeq e_{[SI]}([I])=[I](N-[I])$, that is, the
expected number of ($SI$) pairs is estimated in terms of the number of the expected number of infecteds
$[I]$. This now leads to a self-contained equation in terms of a new approximate variable
$\tilde I$ given by
\begin{equation}
\dot{\tilde{I}}=\tau \tilde I (N-\tilde
I)-\gamma \tilde I.
\label{MFSIS}
\end{equation}
This is the well known compartmental model
for the $SIS$ epidemic.

The same argument can repeated by using a closure at the level of
triples rather then pairs. In this case, the exact pairwise
equations are given by
\begin{eqnarray}
\frac{d[I]}{dt}&=&\tau [SI]-\gamma [I],\\
\frac{d[SI]}{dt}&=&\gamma ([II]-[SI])+\tau
([SSI]-[ISI]-[SI]),\\
\frac{d[II]}{dt}&=&-2\gamma [II]+ 2\tau ([ISI]+[SI]),\\
\frac{d[SS]}{dt}&=&2\gamma [SI]-2\tau [SSI].
\end{eqnarray}
Using the well know closure given by
$[ABC]=\frac{N-2}{N-1}\frac{[AB][BC]}{[B]}$ (see \cite{HouseImproved, Keeling1999}) leads to the following
approximate system
\begin{eqnarray}
\dot {\overline{I}} &=& \tau \overline{SI}- \gamma \overline{I}, \label{TripCloI}\\
\dot {\overline{SI}} &=& \gamma (\overline{II}-\overline{SI}) + \tau (\overline{SSI}-\overline{ISI}-\overline{SI}), \label{TripCloSI}\\
\dot {\overline{II}} &=& -2\gamma \overline{II} + 2\tau
(\overline{ISI}+\overline{SI}) \label{TripCloII},\\
\dot {\overline{SS}} &=& 2\gamma \overline{SI} - 2\tau
\overline{SSI}. \label{TripCloSS}
\end{eqnarray}
The focus now shifts from the derivation of the approximate model to whether and how well
these agree with output from the exact system. More precisely, we will simply consider the difference between the exact solution $[I]$ and the
approximate solutions $\tilde I$ and $\bar I$ in terms of the magnitude of $|[I](t)-\tilde I(t)|$
and $|[I](t)-\bar I(t)|$, and how these depend on population or network size. Numerical
investigation reveals that the difference is of order $1/N$ for closures both at the level of
pairs and triple. This is shown in Fig. 1, where both the exact (Eq. (\ref{eq:Kolmogorov}) with
coefficients given by Eq. (\ref{eq:SIS_coeff_fully_con_graph})) and the approximate systems (i.e.
with pairwise closure given by Eq. (\ref{MFSIS}) and triple closure given by Eqs.
(\ref{TripCloI})--(\ref{TripCloSS})) have been solved numerically. This is somewhat surprising
given that, at least intuitively, it is expected that the closure at the level of triples will be
superior to that at the level of pairs. This numerical result can be made more rigorous at least
for the closure at the level of pairs as shown by
 \cite{BatkaiKissSikolyaSimon, SimonKissIMA}.

\section{Interpreting pairwise equations and closures in terms of moments}

\subsection{Equations for the moments}

In this Section the analysis will focus on the derivation of the
equations for the moments and the interpretation of the pairwise
equations in terms of the moments. Let us define the $j$th
moment associated with the stochastic process as follows
\begin{equation}
y_j(t) = \sum_{k=0} ^N  \left( \frac{k}{N}\right) ^j p_k(t) \quad
\text{or} \quad Y_j(t) = \sum_{k=0} ^N  k^j p_k(t),
\label{DefOfMom}
\end{equation}
where $N^j y_{j}=Y_j$ with $j=1,2, \ldots$. It is straightforward
to derive evolution equations for the moments. For example, the
derivative of the first moment, and in a similar way all others,
can be given in function of other moments upon using the
Kolmogorov equations (\ref{eq:Kolmogorov}). The derivation for the first moment is
outlined below,
\begin{eqnarray*}
\dot{Y}_1(t)&=&\sum_{k=0}^{N}k\dot{p}_{k}=\sum_{k=0}^{N}k(a_{k-1}p_{k-1}-(a_{k}+c_{k})p_k+c_{k+1}p_{k+1})\\
&=&\sum_{k=0}^{N}(ka_{k-1}p_{k-1}-ka_{k}p_k+kc_{k}p_k+kc_{k+1}p_{k+1}).
\end{eqnarray*}
By changing the indices of the summation, plugging in the
corresponding expressions for the transition rates $a_k$ and $c_k$ (Eq. (\ref{eq:SIS_coeff_fully_con_graph}))
and taking into account that $a_{-1}=c_{N+1}=0$ the following
expression holds,
\begin{equation}
\dot{Y}_1(t)=\sum_{i=0}^{N}(\tau (i+i^2)(N-i)-\tau
i^2(n-i)-i^2\gamma+(i^2-i)\gamma)p_i.
\end{equation}
Based on our notations (see, Eq. (\ref{DefOfMom})), the equation above reduces to
\begin{equation}
\dot{Y}_1(t)=\tau N Y_1-\tau Y_2-\gamma
Y_1=(\beta-\gamma)Y_1-\frac{\beta}{N}Y_2,
\label{MomEqY1}
\end{equation}
where $\beta=\tau N$ is the linking relation between mean-field-type and network models. Using a similar procedure, the equation for the second moment
 $Y_2$ can be easily computed and is given by
\begin{equation}
\dot{Y_2} = 2(\beta -\gamma ) Y_2 - 2 \frac{\beta}{N} Y_3 + (\beta + \gamma)Y_1-\frac{\beta}{N}Y_2.
\label{MomEqY2}
\end{equation}
Equations (\ref{MomEqY1}) \& (\ref{MomEqY2}) can be recast in terms of the
density dependent moments $y_j$s to give
\begin{eqnarray}
\dot{y}_1&=&(\beta-\gamma)y_1-y_2,\label{MomEqy1}\\
\dot{y}_2&=& 2(\beta -\gamma ) y_2 - 2 \frac{\beta}{N} y_3 +\frac{1}{N}\left((\beta +
\gamma)y_1-\beta y_2\right).\label{MomEqy2}
\end{eqnarray}
These equations will play a role that is similar to that of the pairwise equations, and similarly,
these are also exact before a closure is applied, at some arbitrarily chosen moment. The above
equations are not closed or self-contained since the second moment depends on the third and an
equation for this is also needed. It is easy to see that this dependence of the moments on higher
moments leads to an infinite but countable number of equations, see
\cite{BatkaiKissSikolyaSimon}. Hence, a closure is needed and in the next Section we will
identify how the classic closure translates to a closure in moments, namely expressing $Y_3$ as a
function of $Y_1$ and $Y_2$.

\subsection{The equivalence between pairwise and moment equations}

To be able to link the pairwise approach to the moment approach it is necessary
to count the expected value of the singles, pairs and triples in terms of the moments.
For a fully connected graph this is straightforward. For example,
based on Eq. (\ref{DefOfMom}), the expected number of infecteds is
given by
\begin{equation}
[I](t)=\sum_{k=0}^{N}kp_{k}(t)=Ny_{1}(t). \label{I_inTermsOfMom}
\end{equation}
Similarly, it is easy to show that similar identities for pairs
and triples can be derived. For example,
\begin{equation}
[SI](t)=\sum_{k=0}^{N}k(N-k)p_{k}(t)=N^2\sum_{k=0}^{N}\frac{k}{N}(1-\frac{k}{N})p_{k}(t)=N^2(y_1(t)-y_2(t)),\label{SI_inTermsOfMom}
\end{equation}
where, $k(N-k)$ simply denotes the number of ($SI$) pairs on a fully connected graph with $N$
nodes and $k$ infected individuals. For triples, the calculations are equally intuitive. For
example, the expected number of $\lbrack SSI \rbrack$ triples can be counted by averaging over
$(N-k)(N-k-1)k$ - the number of $\lbrack SSI \rbrack$ triples in the presence of $k$ infected
nodes. Hence, the following relation holds,
\begin{eqnarray}
\lbrack SSI\rbrack (t)&=&\sum_{k=0}^{N}(N-k)(N-k-1)kp_{k}(t)\notag\\
&=&N^3\sum_{k=0}^{N}\left(1-\frac{k}{N}\right)\left(1-\frac{k}{N}-\frac{1}{N}\right)\frac{k}{N}p_{k}(t)\notag\\
&=&N^3\left(\left(1-\frac{1}{N}\right)y_1(t)+\left(\frac{1}{N}-2\right)y_2(t)+y_3(t)\right).\label{SSI_inTermsOfMom}
\end{eqnarray}
Following the same simple procedure as above, the
following relations hold,
\begin{eqnarray}
[S](t)&=&N(1-y_1(t))\label{S_inTermsOfMom},\\
\lbrack II \rbrack (t)&=&N^2\left(y_2(t)-\frac{1}{N}y_1(t)\right),\label{II_inTermsOfMom}\\
\lbrack SS \rbrack (t)&=&N^2\left(1-\frac{1}{N}+\left(\frac{1}{N}-2\right)y_1(t)+y_2(t)\right),\label{SS_inTermsOfMom}\\
\lbrack ISI \rbrack
(t)&=&N^3\left(-\frac{1}{N}y_1(t)+\left(\frac{1}{N}+1\right)y_2(t)-y_3(t)\right).\label{ISI_inTermsOfMom}
\end{eqnarray}
The results above allows us to test if the equations for the
moments are consistent with the pairwise ones. Starting from the
the equation for the expected number of infecteds and upon using Eq.
(\ref{MomEqy1}), we obtain,
\begin{equation}
[\dot{I}](t)=N\dot{y}_1(t)=N((\beta-\gamma)y_1-y_2)=N(\beta
(y_1-y_2)-\gamma y_1)=\frac{\beta}{N}[SI]-\gamma [I].
\end{equation}
The calculations above can be repeated for all other equations and
these confirm the one-to-one correspondence between the moment and
pairwise equations.

\subsection{Interpreting closures via moments}
The closure at the level of pairs and triples have to equally
translate in a relation between the moments. First, the closure at
the level of pairs is discussed. We look to use the pairwise
closure $[SI]=[S][I]=(N-[I])[I]$ to obtain a relation between the
moments. The simple relation above in terms of the moments
translates to
\begin{equation}
N^2(y_1(t)-y_2(t))=N^2(y_1(t)-y_1^2(t)) \Leftrightarrow
y_2(t)=y_1^2(t).
\end{equation}
For the triple closure the situation is different in that two different triples are closed, albeit
using the same formal relation, and this could potentially lead to two different relations between
the moments. The first triple closure, for the $SSI$ triple, upon using Eqs.
(\ref{SI_inTermsOfMom})--(\ref{S_inTermsOfMom}) \& (\ref{SS_inTermsOfMom}), leads to the following
relation,
\begin{equation}
\left(1-\frac{1}{N}\right)y_1+\left(\frac{1}{N}-2\right)y_2+y_3
 = \frac{N-2}{N-1}\frac{\left(1-\frac{1}{N}+\left(\frac{1}{N}-2\right)y_1+y_2\right)(y_1-y_2)}{1-y_1}.
\end{equation}
The equation above, after some algebra, yields a closure at the level of the moments and $y_3$ can
be given in function of the previous two moments as follows,
\begin{equation}
y_3=-\frac{1}{N} y_1 +(1+\frac{1}{N})y_2 -
\frac{N-2}{N-1}\frac{(y_1-y_2)^2}{1-y_1}.\label{ClassicClosureWithMom}
\end{equation}
It is worth noting that the closure for $[ISI]$ yields the same
closure relation for the third moment.
\section{The new, improved closure and its performance}
The novel closure put forward here is based on the empirical
observation that $p_k(t)$ is well approximated by a binomial
distribution $\mathcal{B}(n,p)$, where $n$ and $p$ depend on time
and will specified in terms of the moments of the distribution.
The first three moments of the binomial distribution can be
specified easily in terms of the two parameters and are as
follows,
\begin{eqnarray}
Y_1&=&np \label{BinMom1}\\
Y_2&=&np+n(n-1)p^2\label{BinMom2}\\
Y_3&=&np+3n(n-1)p^2+n(n-1)(n-2)p^3\label{BinMom3}.
\end{eqnarray}
Using Eqs. (\ref{BinMom1}) \& (\ref{BinMom2}), $n$ and $p$ can be expressed in
term of $Y_1$ and $Y_2$ as follows,
\begin{equation}
p= 1+Y_1- \frac{Y_2}{Y_1}, \qquad  n=\frac{Y_1^2}{Y_1+Y_1^2-Y_2}\label{pAndNinTermOfMom}.
\end{equation}
Plugging the expressions for $p$ and $n$ (Eq. \ref{pAndNinTermOfMom}) into Eq. (\ref{BinMom3}), the closure
for the third moment is found to be
\begin{equation}
Y_3= \frac{2Y_2^2}{Y_1} - Y_2 - Y_1(Y_2-Y_1).
\end{equation}
This relation defines the new triple closure and in terms of the
density dependent moments is equivalent to
\begin{equation}
y_3= \frac{2y_2^2}{y_1} - y_1y_2 + \frac{1}{N} (y_1^2 -y_2).\label{BinomClosure}
\end{equation}
Using the equations for the first moment (\ref{MomEqy1}) the closure at
the level of the pairs yields the following approximate equation
\begin{equation}
\dot x_1 = (\beta -\gamma ) x_1 -  \beta x_1^2.
\end{equation}
Using the equations for the first two moments ((\ref{MomEqy1}) \& (\ref{MomEqy2}))
and the closure at the level of the third moment yields
\begin{eqnarray}
\dot{x_1} &=&  (\beta -\gamma ) x_1 -  \beta x_2, \label{MomEqWithBinomClosure1}\\
\dot{x_2} &=& 2(\beta -\gamma ) x_2 - 2 \beta x_3 + \frac{1}{N}
((\beta + \gamma)x_1-\beta x_2),\label{MomEqWithBinomClosure2}
\end{eqnarray}
where
\begin{equation}
x_3 =\frac{2x_2^2}{x_1} - x_1x_2 + \frac{1}{N} (x_1^2 -x_2).\label{BinomClosureWithX}
\end{equation}
Moreover, we can also define a simplified binomial closure by
neglecting the order $1/N^2$ in the full binomial closure, provided also that $x_1^2-x_2$ is of $\mathcal{O}(1/N)$. This will
lead to
\begin{equation}
x_3 = 3x_1 x_2 - 2x_1^3,
\end{equation}
with all the above in contrast with the the classic triple closure
given by Eq. (\ref{ClassicClosureWithMom}).

The current setup allows us to compare the exact model as given by the Kolmogorov equations
(\ref{eq:Kolmogorov}) with transition rates given by Eq. (\ref{eq:SIS_coeff_fully_con_graph}) to
three different approximate models. The first results from the pairwise closure and is given by
Eq. (\ref{MFSIS}). The second is a direct consequence of the closure at the level of triples and
is given by Eqs. (\ref{TripCloI})--(\ref{TripCloSS}). Finally, the third approximate system
results from the novel binomial closure (see Eq. (\ref{BinomClosure})) with the approximate system
defined by Eqs. (\ref{MomEqWithBinomClosure1})--(\ref{BinomClosureWithX}).  The elegance of this
approach stems from the fact that all numerical results are free from simulations and rely solely
on the time integration of ODE systems. In Fig. 1, plots of the time evolution together with the
approximation errors are given for the classic approximate models corresponding to pairwise and
triple closure. The most significant feature of this plot is the order $1/N$ approximation error independently of the closure. 
The same approximation error suggests qualitative similarity between the two approximate models 
and it somewhat surprising that the closure at the level of triples
produces no immediately obvious benefits or qualitatively different behaviour. In contrast, in
Fig. 2 the performance of the binomial closure is compared to the triple closure.  The
approximation error plot in Fig. 2b is the most significant as it shows numerical evidence that
the binomial closure performs significantly better with error that is of order $1/N^2$. While we
were not yet able to prove analytically this result, for the case of pairwise closure we gave
rigorous proof in \cite{BatkaiKissSikolyaSimon, SimonKissIMA}, where the result is more general 
as apart from estimating the approximate error at the steady state also extends to the approximate error at different times.

The quantitative comparison of the different closure relations is shown in Table 1. Her, the
values of the steady states obtained for the different closures are compared to the exact value of
the stationary prevalence. The advantage of comparing the steady state values is that these can be
given analytically without solving ODEs numerically. The steady state from
the pairwise closure Eq. (\ref{MFSIS}) is
\begin{equation}
\frac{\tilde{I}_{ss}}{N} = 1-\frac{\gamma}{\beta}, \label{ssPC}
\end{equation}
where $\beta=\tau N$. The steady state from the triple closure Eqs.
(\ref{TripCloI})--(\ref{TripCloSS}) is
\begin{equation}
\frac{\overline{I}_{ss}}{N} = (N-1)\frac{(N-2)\beta-N\gamma}{(N-1)(N-2)\beta-N\gamma}.
\label{ssTC}
\end{equation}
 The steady state from the binomial closure Eqs.
(\ref{MomEqWithBinomClosure1})--(\ref{BinomClosureWithX}) is
\begin{equation}
\frac{I^*_{ss}}{N} = \frac{Nq^2-1}{Nq-1}, \label{ssBC}
\end{equation}
where $q=1-\gamma/\beta$. Finally, the value of the steady state from the Kolmogorov equations
(KE) can be derived as
\begin{equation}
\frac{[I]_{ss}}{N} = \frac{\sum\limits_{k=0}^{N-1} (k+1)A_k}{N \sum\limits_{k=0}^{N-1} A_k},
\label{ssKolm}
\end{equation}
where $A_0=1$ and
$$
A_k = \frac{\beta^k (N-1)(N-2)...(N-k)}{\gamma^k (k+1) N^k} , \quad k=1,2,\ldots N-1 .
$$
The approximation errors shown in Table 1 are
$$
\left|\frac{\tilde{I}_{ss}}{N}-\frac{[I]_{ss}}{N}\right|, \quad
\left|\frac{\overline{I}_{ss}}{N}-\frac{[I]_{ss}}{N}\right|, \quad
\left|\frac{I^*_{ss}}{N}-\frac{[I]_{ss}}{N}\right|.
$$
In the Table these values are multiplied by 1000 to get numerical values closer to 1. Using the analytical expressions above it is easy to show that
$\frac{\overline{I}_{ss}}{N} \rightarrow \frac{\tilde{I}_{ss}}{N} $ and $\frac{I^*_{ss}}{N}  \rightarrow \frac{\tilde{I}_{ss}}{N}$ when $N\rightarrow \infty$. However, we did not yet manage to formally prove that
$$
\left|\frac{\tilde{I}_{ss}}{N}-\frac{[I]_{ss}}{N}\right|/N, \quad
\left|\frac{\overline{I}_{ss}}{N}-\frac{[I]_{ss}}{N}\right|/N, \quad
\left|\frac{I^*_{ss}}{N}-\frac{[I]_{ss}}{N}\right|/N^2
$$
all tend to a non-zero number as $N$ gets large and thus analytically confirming the $\mathcal{O}(1/N)$ and $\mathcal{O}(1/N^2)$ errors.

\section{Discussion}

In this paper, we proposed a novel/new closure that produces qualitatively different results when
compared to the classic closure. Namely, the difference between the exact system from the
solution of the approximate model that results from the binomial closure is smaller (i.e.
$\mathcal{O}(1/N^2)$ compared to $\mathcal{O}(1/N)$) and this has been illustrated via numerical examples. The paper has also
identified the link between the equations of the moments and pairwise equations. In particular, it
is worth noting that in terms of moments, the pairwise equations reduce to two simple ODEs
and this is in contrast with four equations for the pairwise model. However, this apparent discrepancy
is not as marked since the pairwise model can be easily reduced to three equations by noting that,
for example, $[SS]=(N-1)(N-[I])-[SI]$. This still leaves one extra equation in the pairwise model
when compared to the equivalent system in terms of moments. This suggests that the three equations
are not completely independent in some sense that is not easy to define or pinpoint.

The present study also allows us to make a few more important observations. For example, one can
attempt to close the moment equations at the second moment but, without using the equivalent
closure that follows from the pairwise model. This can be achieved by considering the limiting case of 
$n = N$ and expressing the second moment in terms of the first. This will lead to a new closure at the level
of pairs but, this performs similarly to the classic closure at the pair level.

The procedure that we illustrate in this paper can be generalised and applied within the
general context of model reduction, where Kolmogorov-type systems are reduced 
to approximate models with far fewer equations.  Our method is based on a semi-heuristic argument
which relies on the assumption that the distribution of the stochastic process in time can be well-approximated 
by some theoretical distribution with time-dependent parameters.  Moreover, the distribution has to be such that there is a unique, one-to-one correspondence between the moments and the parameters of the distribution. If the conditions above hold, the theoretical distribution
will determine the exact type of closure that is necessary to obtain a reduction in the number of equations. The method also allow us to
close at different levels or moments and to explore the benefits of closing at higher order moments. While the initial results and applicability of our method is promising (see example below for the homogeneous random graph case), we are yet to fully explore the applicability of our method or to identify the sufficient and/or necessary requirements that guarantees that the proposed closure or reduction procedure will be applicable. It is worth noting that here the main focus in not on the approximation of simulation results corresponding to a stochastic process but rather the model reduction where the Kolmogorov equations are already given.

Our result about the new closure relation applies directly to homogeneous random graphs provided
that the Kolmogorov equations for this case take the form given in Eq. (KE) with coefficients defined in Eq.
(\ref{eq:SIS_coeff_hom_ran_graph}). The formulation of the result is even more straightforward if
the following slight modification for the transmission rate $a_k$ is introduced
$$
a_k=\tau n k\frac{N-k}{N} .
$$
In this case, the coefficients $a_k$ and $c_k$ in Eq.
(\ref{eq:SIS_coeff_fully_con_graph}) and Eq. (\ref{eq:SIS_coeff_hom_ran_graph}) are formally equivalent, i.e.  both can be
written as $a_k=\frac{\delta}{N} k(N-k)$ with $\delta=\beta$ for the complete graph, and $\delta = \tau n$ for
the homogeneous random graph. It is worth noting that, in the case of the complete
graph, $\tau $ is of order $1/N$ and $\beta$ is independent of $N$, while for the homogeneous
random graph, $\tau$ is independent of $N$. However, in both cases the coefficient $a_k$ can be
given as the product of $k(N-k)$ with a $1/N$ term. Therefore in the case of a homogeneous random
graph the equation for the moments can be easily obtained from Eqs. (\ref{MomEqy1}) \&
(\ref{MomEqy2}) by writing $\tau n$ instead of $\beta$. Thus the ODE system with the binomial closure
for the homogeneous random graph is the same as Eqs.
(\ref{MomEqWithBinomClosure1})--(\ref{BinomClosureWithX}) following a change that again replaces $\beta$ with $\tau n$.

While there are some limitations to the straightforward extension of our results to arbitrary graphs
and/or dynamics, the present paper provides the basis for some more rigorous analysis as well as
the possibility of deriving new closures. As such, this study offers a new or different direction 
that can lead to further progress in the area of approximating stochastic processes independently of whether these involve graphs.

\section*{Acknowledgements} P\'eter L. Simon acknowledges support from OTKA (grant no. 81403).
Funding from the European Union and the European Social Fund is also acknowledged (financial
support to the project under the grant agreement no. T\'AMOP-4.2.1/B-09/1/KMR.).

\newpage
\begin{table}
\begin{tabular}{|l|l|l|l|l|} \hline
N & 100    & 200     & 400 & 800   \\
\hline
1000 $\times$ pair appr. err. & 6.9486 &  3.4008 & 1.6832 & 0.8374 \\
\hline
1000 $\times$ triple appr. err. & 1.2355 &  0.5729 & 0.2763 & 0.1357 \\
\hline
1000 $\times$ binomial appr. err. & 0.1689 &  0.0395 & 0.0096 & 0.0024 \\
\hline
\end{tabular}
\caption{Numerical values for the discrepancy between the exact and approximate solutions for
different network sizes.}
\label{tab:tab}
\end{table}


\newpage
\begin{figure}[ht]
\centering \subfigure[]{
\includegraphics[scale=0.45]{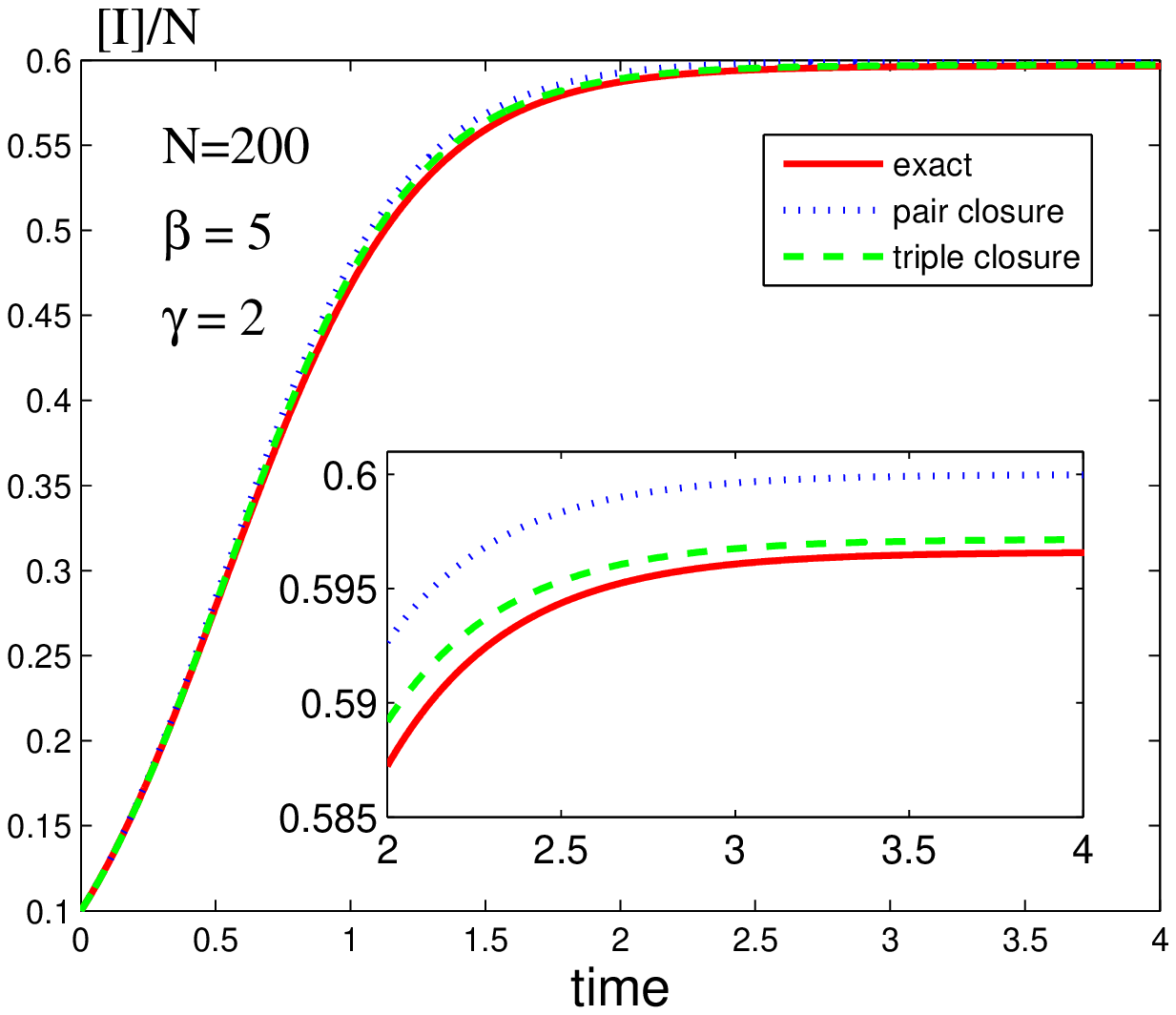}
\label{fig:fig1a} } \subfigure[]{
\includegraphics[scale=0.45]{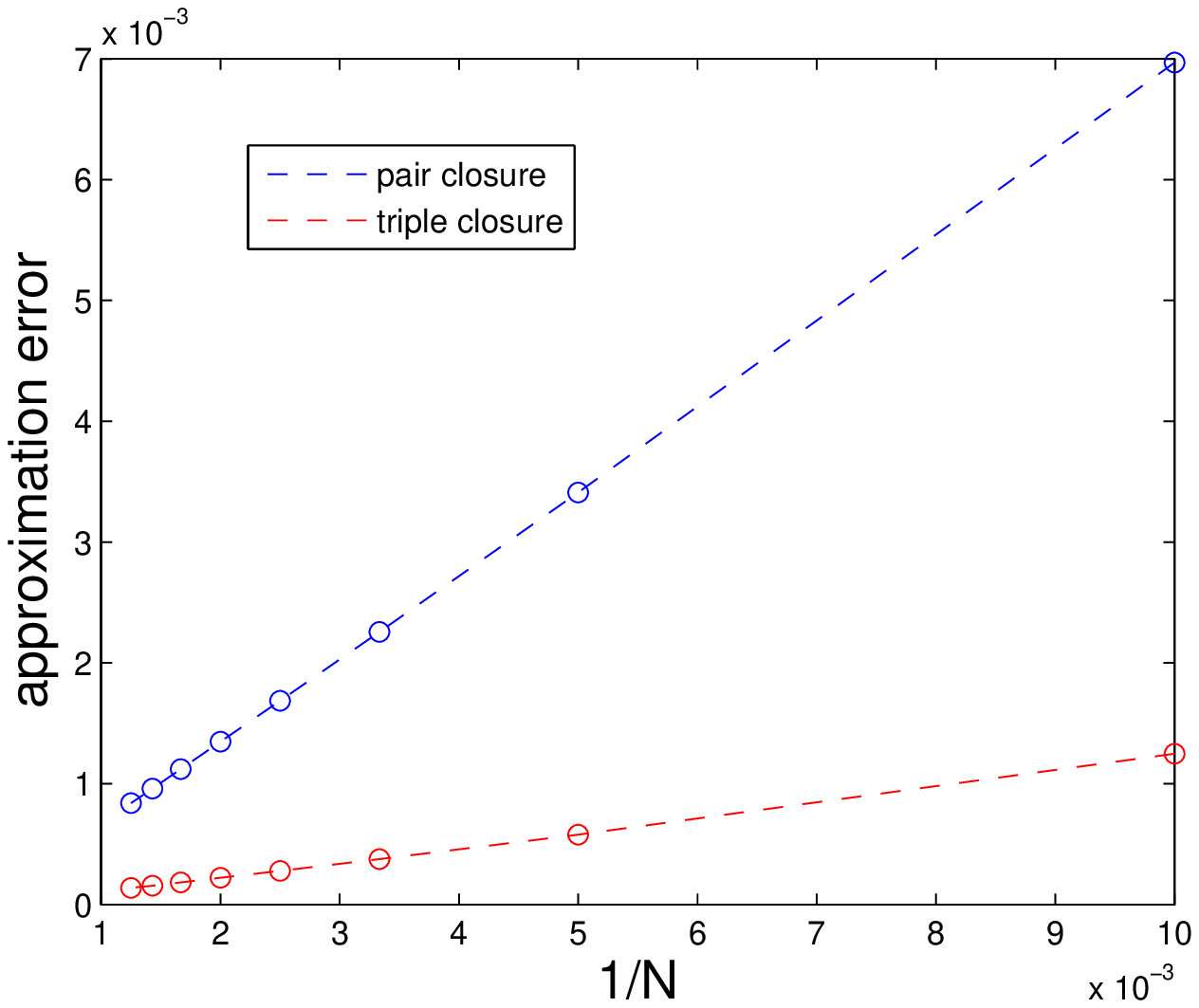}
\label{fig:fig1b} } \caption[Optional caption for list of figures]{(a) Time evolution of the
fraction infected (prevalence) based on the exact model and two different closures; pair and
triple level. Here, $N=200$, $\gamma=2$ and $\beta=5$. (b) Approximation error (absolute value of
the difference between the exact and two approximate models at steady state) plotted for different
system sizes for the same transmission and recovery parameter values as at (a). } \label{fig:fig1}
\end{figure}

\newpage
\begin{figure}[ht]
\centering \subfigure[]{
\includegraphics[scale=0.45]{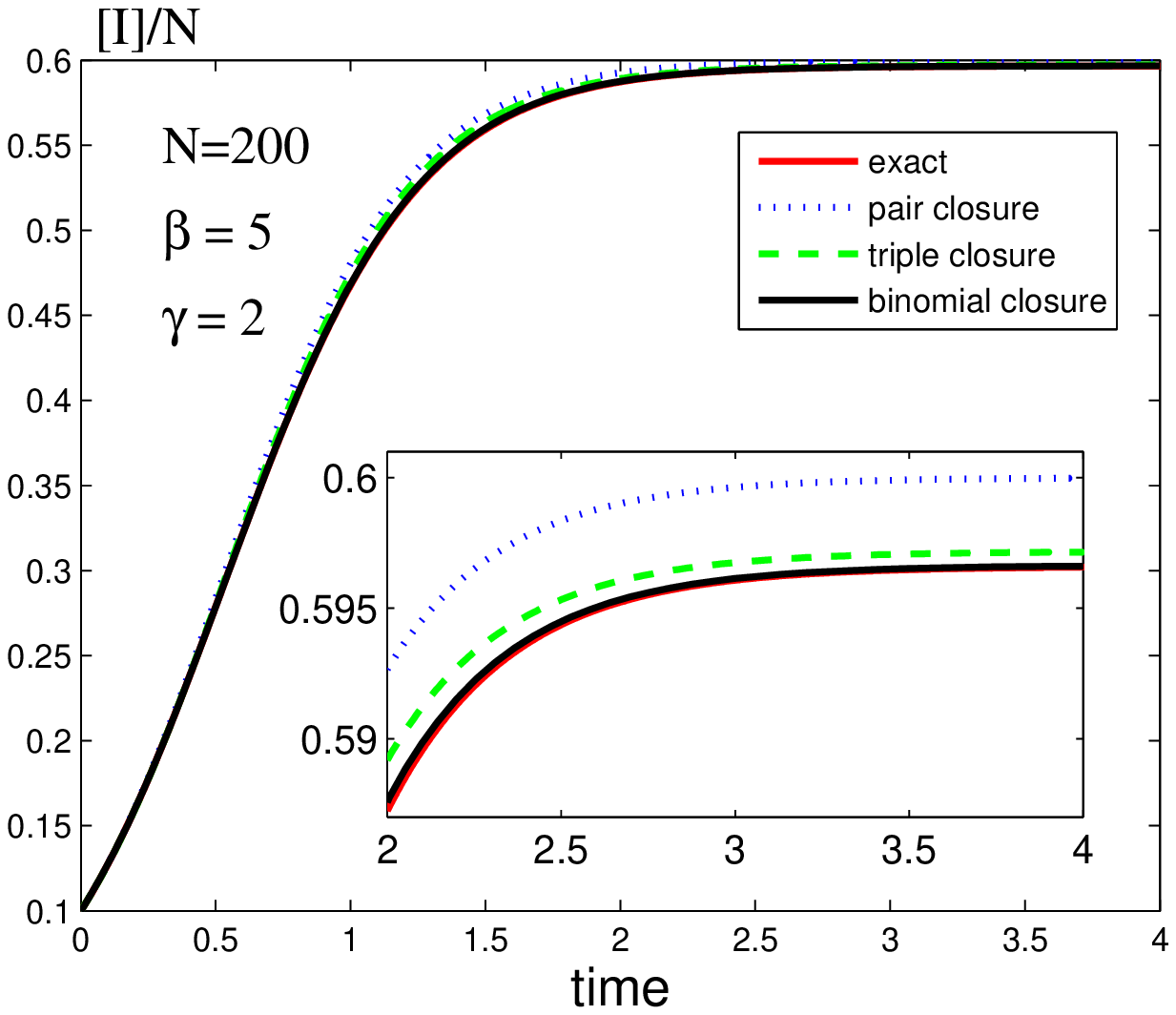}
\label{fig:fig2a} } \subfigure[]{
\includegraphics[scale=0.45]{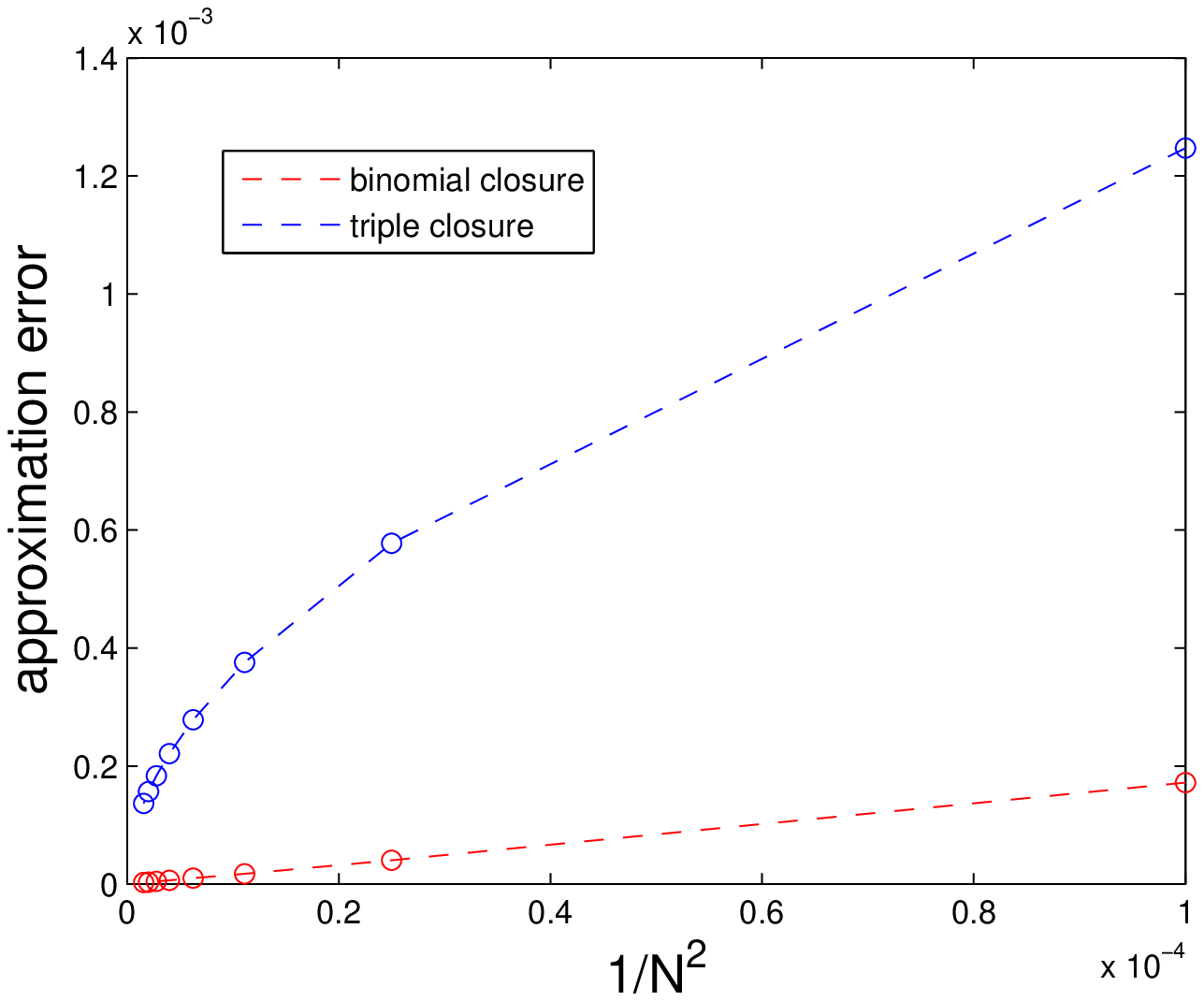}
\label{fig:fig2b} } \caption[Optional caption for list of figures]{(a) Time evolution of the
fraction infected (prevalence) based on the exact model and three different closures; pair, triple
level and the binomial closure. Here, $N=200$, $\gamma=2$ and $\beta=5$. (b) Approximation error
(absolute value of the difference between the exact and approximate models at steady state given
by the triple and binomial closures) plotted for different system sizes for the same transmission
and recovery parameter values as at (a). It is worth noting that the $x$-axis in (b) is in terms of $1/N^2$ as opposed to $1/N$ as in Fig. 1(b).}
\label{fig:fig2}
\end{figure}

\end{document}